\documentclass[10pt,fleqn]{article}

\usepackage{amsfonts,amssymb,amsmath}
\usepackage{graphicx} 
\usepackage{color}
\usepackage[bookmarksnumbered, plainpages, backref, colorlinks]{hyperref}

\newcommand{\R}{\mathbb R}

\newcommand{\ds}{\displaystyle}
\newcommand{\la}{\langle}
\newcommand{\ra}{\rangle}

\newcommand{\qed}{\hfill $\Box$}

 \def\reff#1{\mbox{\rm(\ref{#1})}}

\newtheorem{theorem}{Theorem}
\newtheorem{lemma}{Lemma}
\newtheorem{proposition}{Proposition}

\newtheorem{remark}{Remark}
\newtheorem{example}{Example}
\newtheorem{corollary}{Corollary}

\newcommand{\GP}{{\rm I\hspace{-.15em}P}}

 \newcommand{\xdba} [2]
    {\left\langle {#1}\,,\,{#2} \right\rangle   }

\def\div{\mbox{\rm div\hspace{1mm}}}
\def\dist{\mbox{\rm dist\hspace{1mm}}}

\def\norm#1#2{\Vert\,#1\,\Vert_{#2}}
\def\hnorm#1#2{\vert\,#1\,\vert_{#2}}

\title{Coupling of finite element and boundary element methods with regularization for a nonlinear interface problem  
with nonmonotone set-valued transmission conditions}
\author{J. Gwinner \& N. Ovcharova \\ 
 Department of Aerospace Engineering \\
 Universit\"{a}t der Bundeswehr M\"{u}nchen \\ 
Neubiberg - Munich, Germany} 

\begin{document}
\date {} 
\maketitle

\begin{center} 
{\it Dedicated to  %our dear colleague 
Professor E.P. Stephan \\
on the occasion of his 75th birthday}
\end{center}

\begin{abstract}
For the first time, a nonlinear interface problem on an unbounded domain  
with nonmonotone set-valued transmission conditions is analyzed.
The investigated problem involves a nonlinear monotone partial differential 
equation in the interior domain and the Laplacian in the exterior domain.
Such a scalar interface problem models nonmonotone frictional contact of 
elastic infinite media. The variational formulation of the interface problem leads to a hemivariational inequality, which lives on the unbounded domain,
and so cannot be treated numerically in a direct way.
By boundary integral methods the problem is transformed and a novel 
hemivariational inequality (HVI) is obtained that lives on the interior domain and on the coupling boundary, only. Thus for discretization the coupling of 
finite elements and boundary elements is the method of choice.
In addition smoothing techniques of nondifferentiable optimization are adapted and the nonsmooth part in the HVI is regularized. Thus we reduce the original variational problem to a finite dimensional problem that can be solved by standard optimization tools. We establish not only convergence results for the total approximation procedure, but also  an asymptotic error estimate for the regularized HVI.
\end{abstract}

\\[1ex]

\section{Introduction}

This paper presents a novel coupling of finite element  and boundary element
methods combined with a regularization procedure for the solution of a 
nonlinear scalar interface problem on an unbounded domain  
with nonmonotone set-valued transmission conditions
that models nonlinear contact problems with nonmonotone friction in infinite elastic media. Such contact problems arise in various fields of science and technology; let us mention geophysics, see e.g. \cite{scholz2019}, soil mechanics,
in particular soil-structure interaction problems, see e.g. \cite{GhIo2009},
 and civil engineering of underground structures, see e.g. \cite{wang2018soil}.

In this paper we employ various mathematical techniques for the solution of the interface problem that consists of a nonlinear monotone partial differential equation in a bounded domain and the Laplace operator in the exterior domain coupled by set-valued nonmonotone transmission conditions. 
Using  Clarke generalized differentiation \cite{Clarke} we describe this coupled boundary value problem as a nonlinear 
{\it hemivariational inequality} (HVI). Further by singular boundary integral methods \cite{HsWe-2008} we reduce this problem to a HVI that lives on the bounded domain and the coupling boundary, only. Thus for discretization we can develop a 
{\it BEM/FEM coupling method}. 
In addition we adapt from \cite{OvGw-2014} smoothing techniques of 
{\it nondifferentiable optimization} and regularize the nonsmooth part in the HVI to arrive at a finite dimensional problem that can be solved by standard optimization tools. We do not only provide convergence results for the total approximation procedure, but also show an error estimate for the regularized HVI that %compares to
gives the same convergence order as the error estimate of BEM/FEM coupling 
for the monotone Signorini contact problem in \cite{CCJG-1997}.

The {\it coupling of finite element and boundary element methods} combining the best of both "worlds" \cite{Zienk-1979} provides nowadays a very effective tool for the numerical solution of boundary value problems in physics and engineering. Indeed, the boundary element method is better suited to problems in which the
domain extends to infinity but is usually confined to regions in which the governing equations are linear and homogeneous. On the other hand, the finite element method is restricted to problems in bounded domains but is applicable to problems in which the material properties are not necessarily homogeneous and nonlinearity may occur. 
This method was originally proposed by engineers (see e. g. Zienkiewicz et al. \cite{Zienk-1979}). The mathematical analysis goes back to Brezzi,
Johnson, and N\'{e}d\'{e}lec \cite{BrJoNe1978,BrJo1979,JoNe1980}
 and was extended by Wendland \cite{Wend-1986,wend-1988}. The 
{\it symmetric coupling} of finite element and boundary element methods
is due to Costabel \cite{Cost-1987,Cost-1988}. 
For the coupling of finite element and boundary element methods for various
{\it nonlinear  interface problems} we point out the papers \cite{CoSt-1990,Ste92,gatih92a,gatih95,MaiSte-2005,CGMS-2006,GMSS-2011,AFFKMP-2013,GMS-2021} and for a comprehensive exposition of finite element and boundary 
coupling refer to \cite[Chapter 12]{GwiSte-2018}.

The theory of {\it hemivariational inequalities} has been introduced and studied since 1980s by Panagiotopoulos \cite{Pan-1993}, as a generalization of variational inequalities with an aim to model many problems coming from mechanics when the energy functionals are nonconvex, but locally Lipschitz,
so the Clarke generalized differentiation calculus \cite{Clarke} 
can used, see \cite{Naniewicz,GMDR-2003,GoeMot-2003}.
For  more recent monographs on hemivariational inequalities with  application to contact problems we refer to 
\cite{Ochal,SofMig-2018}. 
In parallel with the mathematical analysis of hemivariational
inequalities the interest in efficient and reliable numerical methods for their solution constantly increases. The classical book on the finite element method for hemivariational inequalities
is the monograph of Haslinger et al. \cite{HaMiPa-1999}.
More recent work on numerical solution of HVIs modeling 
contact problems with nonmonotone friction, adhesion, cohesive cracks, and delamination for elastic bodies on bounded domains is contained in the papers
\cite{BHM,Kov-11,HK2,LPS2,Nes-12,HaSoBa-2017,HanSof-2019};
see also \cite{Wohl-2011} for variationally consistent discretization schemes and numerical algorithms for contact problems.

The plan of the paper is as follows. 
The next section 2 collects some basic notions of Clarke's generalized differential calculus that are needed for the analysis of the nonmonotone
transmission conditions. Then we descibe the interface problem under study. 
We settle the issues of existence and uniqueness in a general functional analytic framework.
Section 3 provides a first equivalent  weak variational formulation  of the interface problem  
in terms of a hemivariational inequality (HVI). Since this 
HVI lives on the unbounded domain 
(as the originall problem), this hemivariational formulation cannot numerically 
treated directly and therefore provides only 
an intermediate step in the solution procedure.
Section 4 employs boundary integral analysis to transform the interface 
problem to a HVI that lives on the interior bounded domain and the interface
boundary, only, and so is amenable to numerical treatment. 
Section 5 uses regularization techniques of nonsmooth optimization and  
derives a regularized version of the HVI that becomes a variational equality. 
Section 6 turns to numerical analysis of the regularized HVI, employs the Galerkin boundary element/finite element method,  and presents an 
asymptotic error estimate.
The final section 7 shortly summarizes our findings, gives some concluding remarks,  and sketches some directions of further research.

\section {Some preliminaries and the interface problem} \label{inter}

Let us first recall the central notions of Clarke's generalized differential calculus \cite{Clarke}, before we pose our interface problem.   
Let $X$ be a (real) Banach space, let $f: X \to \R$
be a locally Lipschitz function. Then  
$$ f^0(x;z) := \ds \lim \sup_{ y \to x; t \downarrow 0}
  \frac{f(y+tz) - f(y)}{t} ~ ~ x, z  \in X,  $$
is called the {\it generalized directional derivative} of $f$ in the direction $z$. 
Note that the function $z \in X \mapsto f^0(x;z)$ is finite, positively homogeneous, and sublinear, hence convex and continuous; further, the function 
$(x,z) \mapsto f^0(x;z)$ is upper semicontinuous.
The {\it generalized gradient} of the function $f$ at $x$,
 denoted by (simply) $\partial f(x)$,
is the unique nonempty $\text{weak}^*$ compact convex subset of
the dual space $X'$, whose support function is $f^0(x;.)$.
Thus 
\begin{eqnarray*} 
&& \xi \in \partial f(x) \Leftrightarrow f^0(x;z) \ge \la \xi,z \ra, \, \forall z \in X , \\ 
&& f^0(x;z) = \max \{ \la \xi,z \ra ~:~ \xi \in \partial f(x) \}, \, \forall z \in X \,. 
\end{eqnarray*}
When $X$ is finite dimensional, then, according to Rademacher's theorem, $f$ is differentiable almost everywhere, and the  generalized gradient of  $f$ at a point $x\in \R^n$ can be characterized by
$$
\partial f(x)= \textnormal{co} \,
 \{\xi\in \R^{ n} \, : \, \xi = \displaystyle \lim _{k \to \infty}
\nabla f(x_k), \; x_k \to x, \,f \, \textnormal{is differentiable at }\, x_k\},
$$
where "co" denotes  the convex hull.\\[2ex]
In this paper we treat the following interface problem. 
Let $\Omega \subset \R^d ~(d \ge 2)$ be a bounded domain with Lipschitz
 boundary $\Gamma$. To describe mixed transmission conditions, we  
assume that the boundary $\Gamma$ splits into two non-empty, open  disjoint parts $ \Gamma_s$ and $\Gamma_t$ such that $\Gamma = 
\textnormal{cl  } \Gamma_s \cup \textnormal{cl } \Gamma_t$.
Let $n$ denote the unit normal on $\Gamma$
defined almost everywhere
pointing from $\Omega $ into $\Omega^c :=\R^d\setminus\overline{\Omega} $.

In the interior part $\Omega$, we consider the nonlinear partial differential equation
\begin{equation}\label{a1}
\div \Bigl( p(| \nabla u |)\cdot \nabla u \Bigr) + f_0 = 0 
\qquad \mbox{in } \; \Omega,
\end{equation}
where $p:[0,\infty) \to [0,\infty)$ is a continuous
function with $ t \cdot p(t) $
being  monotonously increasing with $t$.

In the exterior part  $\Omega^c$, we consider the Laplace equation
\begin{equation}\label{a2}
\Delta u=0 \qquad \mbox{in }\; \Omega^c
\end{equation}
with the radiation condition at infinity $( |x|\to\infty) $

\begin{eqnarray}
\label{a3}     %\nonumber
u(x) =  
    \left\{ \begin{array}{ll}
     a + o(1)  & \mbox{if } d=2 \,, \\[0.5ex]  
      O(|x|^{2-d}  )  & \mbox{if } d > 2 \,, 
    \end{array} \right.
  \end{eqnarray}

where $a$ is a real constant for any $u$, but may vary with $u$.
Let us write $u_1:=u|_\Omega $ and  $u_2:=u|_{\Omega^c} $, then the tractions on the coupling boundary $\Gamma$ are given by the traces of
$ p(|\nabla u_1|) \frac{\partial u_1}{\partial n}  $ and
$-\frac{\partial u_2}{\partial n} $, respectively. 
 
We prescribe classical transmission conditions on $\Gamma_t$,
\begin{equation}\label{a4}
u_1|_{\Gamma_t} = u_2|_{\Gamma_t}+ u_0|_{\Gamma_t} \quad\mbox{and}\quad
p(|\nabla u_1|)\left.\frac{\partial u_1}{\partial n}\right \vert_{\Gamma_t}
=  \left.\frac{\partial u_2}{\partial n}\right\vert_{\Gamma_t} + t_0|_{\Gamma_t} ,
\end{equation}
and  on $\Gamma_s$, analogously for the tractions, 
\begin{equation} \label{a5_1}
p(|\nabla u_1|)\left.\frac{\partial u_1}{\partial n}\right\vert_{\Gamma_s}
=  \left.\frac{\partial u_2}{\partial n}\right\vert_{\Gamma_s} + t_0|_{\Gamma_s}
\end{equation}
and the generally nonmonotone, set-valued transmission condition,   
\begin{equation}  \label{a5}
 p(|\nabla u_1|) \left. \frac{\partial u_1}{\partial n}\right\vert_{\Gamma_s}  
  \in \left. \partial j(\cdot, u_0  + (u_2 - u_1)\right\vert_{\Gamma_s}) 
 \,.
\end{equation}

Here  the function $j: \Gamma_s \times \R \to \R $ is such that $j(\cdot, \xi):\Gamma_s\to \R$ is  measurable on $\Gamma_s$ for all $\xi \in \R$ and $j(s, \cdot) : \R \to \R$ is locally Lipschitz  for almost all (a.a.) $s\in \Gamma_s$  with  
 $\partial j(s,\xi):=\partial j(s,\cdot) (\xi)$, the generalized gradient of $j(s,\cdot)$ at $\xi$.

Further, we require the following growth condition on the 
so-called superpotential $j$:
There exist positive constants $c_3$ and $c_4$ such that for a.a. $s \in \Gamma_s$, all $\xi \in \R$ and for all 
$\eta \in \partial j(s, \xi)$ the following inequalities hold
\begin{equation} \label{as-j} % (H(j)) 
          (i)\quad |\eta| \leq c_3 (1+|\xi|) \,, 
          (ii)\quad \eta ~\xi \geq - c_4|\xi| \,.
     \end{equation}

 \begin{remark}
 The friction condition of the paper includes Tresca friction or following
 \cite{MaiSte-2005}  Coulomb friction with given nonnegative friction force $g$.
  Indeed, choose $j(\cdot,\xi) = g(\cdot) |\xi|$ for $d=2$, then 
  \begin{eqnarray*}
  \partial j (\cdot, \xi) =
     \left\{ \begin{array} {lll}
- g & \mbox{ if } & \xi < 0 \\[0.5ex]
[-g, g]  & \mbox{ if } &  \xi = 0 \\[0.5ex]
g  & \mbox{ if} &  \xi > 0 
\end{array} \right.
\end{eqnarray*}
is monotone set-valued and \reff{a5} becomes
\[
| p(|\nabla u_1|)  \frac{\partial u_1}{\partial n}|  
\le  g (u_0 + (u_2 - u_1)\vert_{\Gamma_s})   \mbox{ on } \Gamma_s \,.
\]
  \end{remark}

Given data $f_0\in L^2(\Omega)$, $u_0\in H^{1/2}(\Gamma)$, and
$t_0 \in L^2(\Gamma)$ together with
\begin{equation} \label{as1}
\int_\Omega f_0 \, dx + \int_\Gamma t_0 \, ds = 0, \; \mbox{if} \; d=2, 
\end{equation}
we are looking for  $u_1\in H^1(\Omega)$ and $u_2\in H^1_{loc}(\Omega^c )$ satisfying \reff{a1}--\reff{a5}
in a weak form.\\[2ex]
%%%%% schwache Form vorweg:
To conclude this section, we now discuss the existence and uniqueness of a weak solution to this interface problem in functional analytic terms.
To this end, let
$X:= L^2 (\Gamma_s)$ and introduce the real-valued locally Lipschitz functional
$$
J(y) := \int_{\Gamma_s} j(s,y(s)) ~ ds\,,  \qquad y \in X \,.
$$
Then by Lebesgue's theorem of majorized convergence,
$$
J^0(y;z)= \int_{\Gamma_s} j^0(s,y(s);z(s)) ~ ds\,, 
 \qquad (y,z) \in X\times X  \,,
$$
 where $j^0(s, \cdot~;~\cdot)$ denotes the generalized directional derivative of $j(s,\cdot)$.

As we shall see in the subsequent sections, the weak formulation of the problem  \reff{a1}--\reff{a5} leads, in an abstract setting, to a  hemivariational inequality (HVI) with a nonlinear operator  ${\cal A}$ and the nonsmooth functional $J$, namely, we are looking for some  
$\hat v\in {\cal C} $ such that
\begin{equation}\label{a6}
{\cal A}(\hat v) (v- \hat v) 
+ J^0(\gamma \hat v; \gamma v- \gamma \hat v)  %\varphi(\hat{v},v)
\ge \lambda(v-\hat v)
\qquad \forall v \in {\cal C}.
\end{equation}
Here ${\cal C}\ne \emptyset $ is a closed convex subset of a Banach space
$E$, the nonlinear operator ${\cal A} : E\to E^* $ is a monotone operator,
and $\gamma := \gamma_{E \to X}$ denotes the linear continuous trace operator,
and the linear form $\lambda$ belongs to the dual $E^*$.
Similar to \cite{CCJG-1997}, the
operator ${\cal A}$ consists of  a nonlinear  monotone differential operator (as made precise below) that results from the PDE \reff{a1} in the bounded domain $\Omega$ and the  positive definite Poincar\'e--Steklov operator on the boundary $\Gamma$ of $\Omega$ that stems from 
 the exterior problem \reff{a2}-\reff{a3} and can be represented by the boundary integral operators of potential theory. Thus it results that the operator ${\cal A}$ is strongly monotone with some monotonicity 
constant $c_{\cal A} > 0$ and Lipschitz continuous on bounded sets. On the other hand, by the compactness of the trace map $\gamma$, the  
real-valued upper semicontinuous bifunction 
$$\psi(v,w) :=  J^0(\gamma v; \gamma w- \gamma v)
  \,, \forall (v,w) \in E\times E  $$
can be seen to be pseudo-monotone, see 
\cite[Lemma 1]{OvGw-Rassias-14}, \cite[Lemma 4.1]{GwOv-2015}.
The latter result also shows a linear growth of $\psi(\cdot,0)$.
This and the strong monotonicity of ${\cal A}$ imply coercivity. 
Therefore by the theory of pseudo-monotone VIs
\cite[Theorem 3]{Gwi-81}, \cite{Zei-2}, see \cite{GwOv-2015} 
for the application to HVIs, we have solvability  of  \reff{a6}.

Further suppose that the generalized directional derivative $J^0$
satisfies the  one-sided Lipschitz condition: There exists $c_J > 0$
such that 
\begin{equation}\label{cc1}
J^0(y_1; y_2 - y_1) + J^0(y_2; y_1 - y_2)
\le c_J \| y_1 - y_2 \|_X^2 
\qquad \forall y_1,y_2 \in X \,.
\end{equation}
Then the smallness condition
\begin{equation}\label{cc2}
c_J \| \gamma \|^2_{E \to X} < c_{\cal A}
\end{equation}
 implies unique solvability of  \reff{a6}, see e.g. 
\cite[Theorem 5.1]{Ov2017} and \cite[Theorem 83]{SofMig-2018}.

It is noteworthy that under the smallness condition  \reff{cc2}
together with  \reff{cc1}, fixed point arguments 
\cite{Capa-2014} or the theory of set-valued pseudo-monotone operators \cite{SofMig-2018} are not needed, but simpler monotonicity arguments are sufficient to conclude unique solvability. Thus the compactness of the trace map $\gamma$ is not needed either.
In fact, \reff{a6} can be framed as a
{\it monotone equilibrium problem} in the sense of Blum-Oettli \cite{BlOe94}:

\begin{proposition} \label{prop1}
 Suppose \reff{cc1} and \reff{cc2}. Then the bifunction
 $\varphi: {\cal C} \times {\cal C} \to \R$ 
defined by 
$$\varphi(v,w) := 
{\cal A}(v) (w- v) 
+ J^0(\gamma v; \gamma w - \gamma v) - \lambda(w-v) $$
has the following properties: \\
$\varphi(v,v) = 0 \,\, \forall v \in {\cal C}$;\\
$\varphi(v,\cdot)$ is convex and lower semicontinuous 
$\forall v \in {\cal C}$; \\
there exists some $\mu > 0$ such that 
$\varphi(v,w) + \varphi(w,v) \le  
- \mu  \|v - w \|^2_E \,\,  \forall v,w \in {\cal C}$ (strong monotonicity).
Moreover, the function $t \in [0,1] \mapsto \varphi(tw + (1-t)v,w)$ is upper semicontinous at $t= 0$ for all $v,w \in {\cal C}$ (hemicontinuity).
\end{proposition}
{\bf Proof.}
Obviously $\varphi$ vanishes on the diagonal and is convex and lower semicontinuous with respect to the second variable.
To show strong monotonicity, estimate
\begin{eqnarray*}
&& \varphi(v,w) + \varphi(w,v) \\
&=& ({\cal A}(v) - {\cal A}(w)) (w- v) \\ 
&& + \, J^0(\gamma v; \gamma w - \gamma v) + J^0(\gamma w; \gamma v - \gamma w) \\ 
&\le & - c_{\cal A} \, \|v - w \|^2_E + 
c_J \| \gamma v - \gamma w \|_X^2  \\
&\le & - (c_{\cal A} - c_J  \| \gamma \|^2_{E \to X} ) \, \|v - w \|^2_E \,. 
\end{eqnarray*}
To show hemicontinuity, it is enough to consider the 
bifunction $ (y,z) \in X \times X \mapsto J^0(y; z - y)$. Then for 
 $ (y,z) \in X \times X $ fixed, $t \in [0,1]$ one has
$$
J^0(y + t(z-y); z - (y + t(z-y))) = (1-t) J^0(y + t(z-y);z-y)
$$
and thus hemicontinuity follows from upper semicontinuity of $J^0$,
$$
\limsup_{t \downarrow 0}  J^0(y + t(z-y);z-y) \le  J^0(y;z-y) \,.
$$    
\qed

Since strong monotonicity implies coercivity and uniqueness,  the fundamental existence result \cite[Theorem 1]{BlOe94} applies to the HVI (\ref{a6}) to conclude the following

\begin{corollary} \label{cor-1}
 Suppose \reff{cc1} and \reff{cc2}.
Then the  HVI (\ref{a6}) is uniquely solvable.
\end{corollary}

Let us comment on the crucial assumptions \reff{cc1} and \reff{cc2}.
The smallness condition  \reff{cc2} says that the strong monotonicity of the 
operator ${\cal A}$ dominates the nonmonotone term $J^0$. 
For a  better understanding of the  one-sided Lipschitz condition  \reff{cc1} 
we insert some of the explanations from  \cite{Ov-Banz} and present  a class of locally Lipschitz functions that satisfies \reff{cc1}.
\\
Let $ f:R \to R$ be  a function such that
\begin{equation}
( \xi^*-\eta^*)\,  (\xi- \eta)  \geq -c_J |\xi -\eta |^2 \quad 
\forall \,  \xi^*\in \partial j(\xi), \; \forall \,  \eta^*\in \partial j(\eta) \label{eq1}
\end{equation}
for any $\xi, \eta \in R$ and  some  $c_J \geq 0$.
From the definition of the Clarke generalized derivative we get
$$
j^0(\xi;\eta-\xi)= \max _{\xi^*\in \partial j(\xi)}  \xi^* \,  (\eta -\xi). 
$$
Rewriting (\ref{eq1}) as
$$
\xi^* \, (\eta- \xi)  +  \eta^* \,  (\xi- \eta) \leq c_J |\xi -\eta |^2
$$
we find 
\[
j^0(\xi;\eta-\xi)+ j^0(\eta;\xi-\eta) \leq c_J|\xi-\eta|^2.
\]
Hence,
\begin{align}
 \psi(v,w)+\psi(w,v) &= 
\int_{\Gamma_s} j^0( \gamma v; \gamma w - \gamma v )\, ds + \int_{\Gamma_s} j^0( \gamma w;  \gamma v -\gamma w )\, ds 
\nonumber \\
& \leq c_J \| \gamma v -\gamma w \|^2_{L^2(\Gamma_s)} \leq c_J \gamma_0 \|v-w\|^2_{E} \label{eq2}
\end{align}
where $\gamma_0$ is the norm of the trace mapping from $E:= H^1(\Omega)$ into $L^2(\Gamma_s)$, and which w. l. o. g. can be considered to be $1$. That's why we simply assume that $c_J$ is sufficiently small, i.e $c_J < c_{\mathcal{A}}$.

Next, we show that if $\partial j$ includes only non-negative jumps, then the condition (\ref{eq1}) is globally satisfied. Whereas for negative jumps, the condition (\ref{eq1}) holds only locally and for details we refer to [39].
\\
\begin{figure}[t] 
\centering
\includegraphics[trim=2cm 10cm 3cm 10cm,clip,width=0.5\textwidth]
{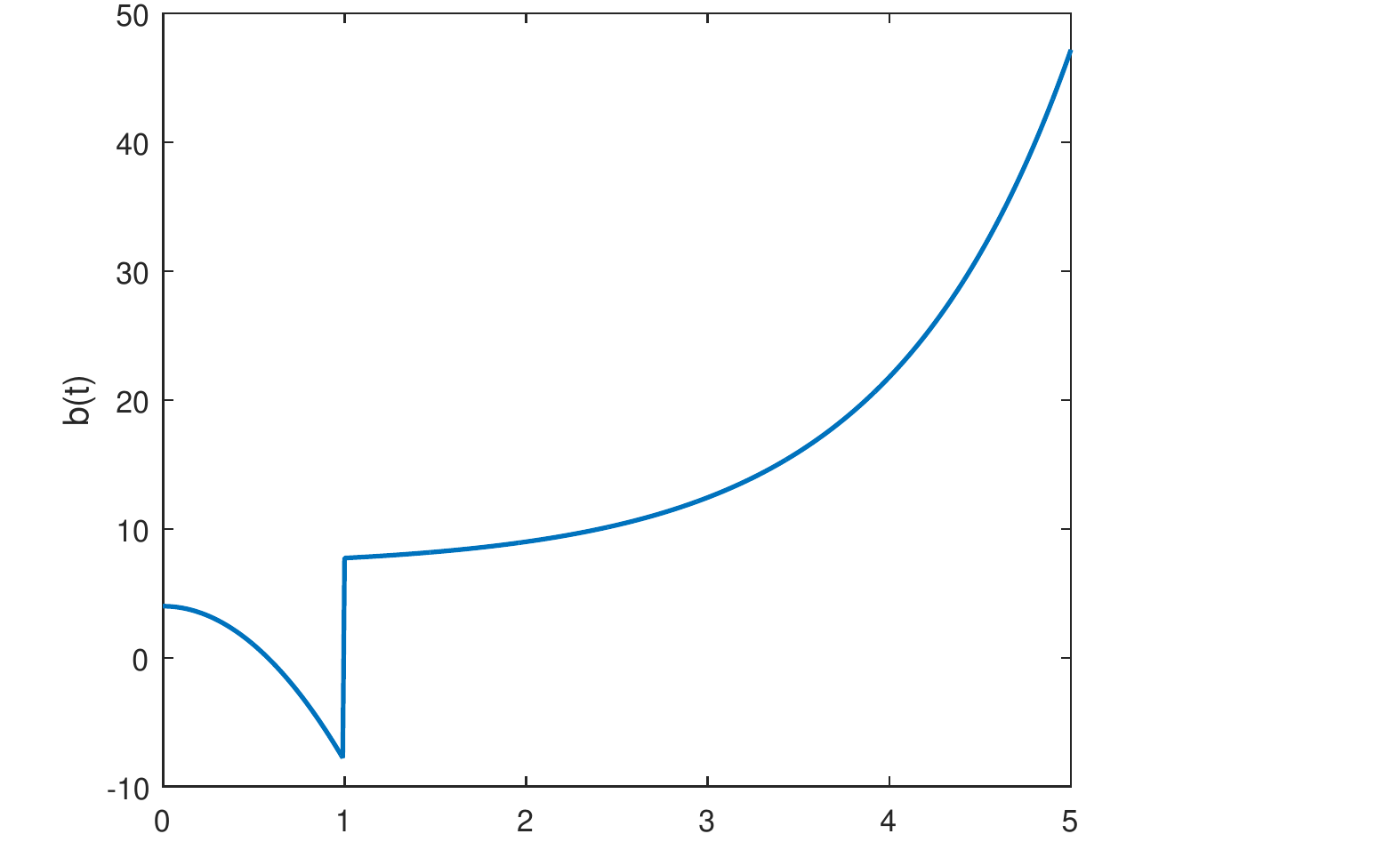} \hspace{-1.5cm}
\includegraphics[trim=2cm 10cm 3cm 10cm,clip,width=0.5\textwidth]{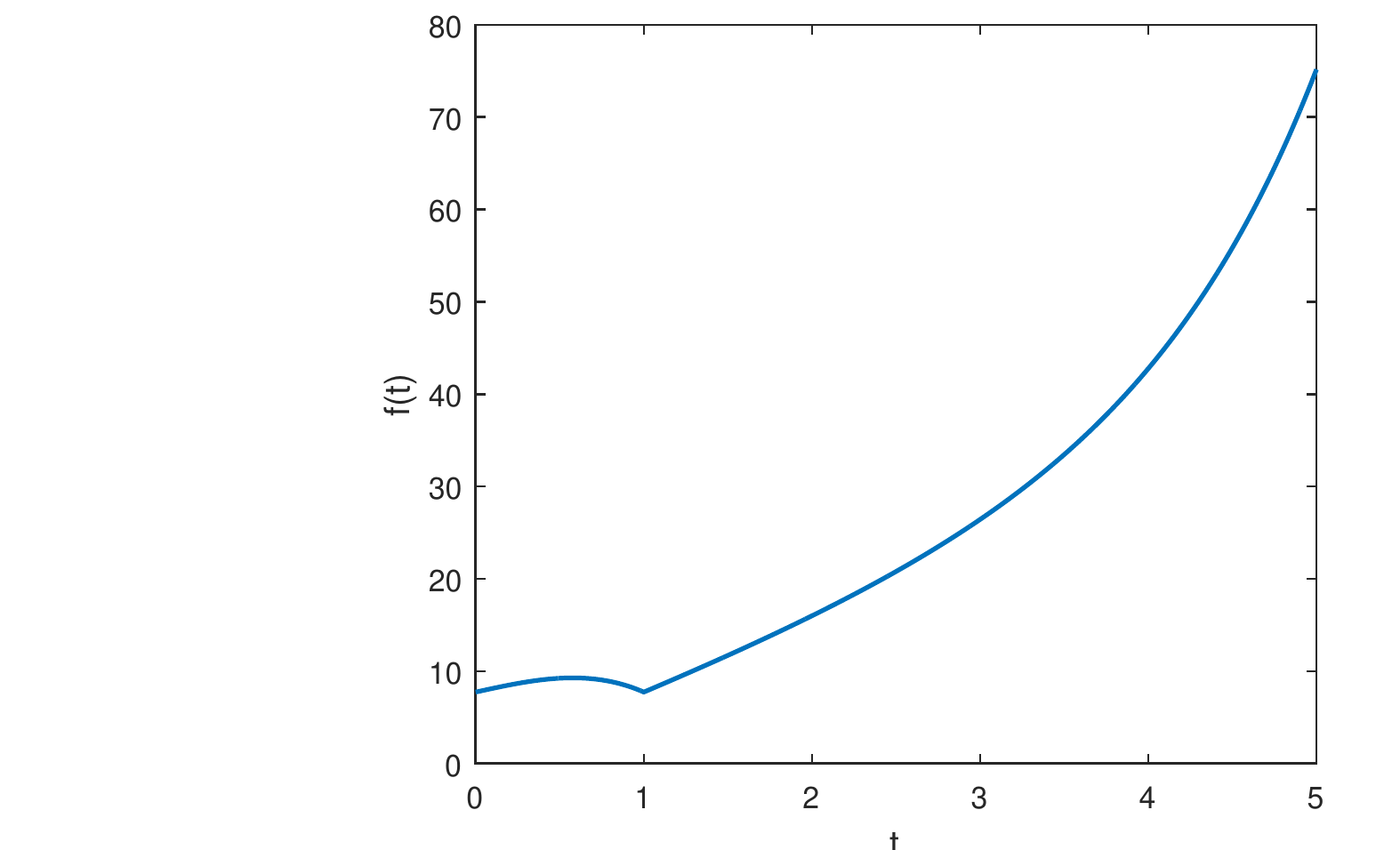}
\caption{An exemplary function $b$ with one positive jump and its anti-derivative  $f(x)= \int_0^x b(t)\,dt$ }
 \label{nonnegative_jump}
\end{figure}
\begin{example}
Let $f :R \to R$, $f(x)=\int_0^x {b(t)}\, dt$, be the anti-derivative of a piecewise Lipschitz function $b$ with finite non-negative jumps at the points $t_i^J$ of discontinuity (see Fig. \ref{nonnegative_jump}).
This means that
$\underline{b}(t_i^J)\leq \overline{b}(t_i^J)$, where  
$$
\underline{b}(t_i^J):=f'(t_i^J-0)= \lim _{h\to 0^-}\frac{f(t_i^J+h)-f(t_i^J)}{h}
$$
and 
$$
\overline{b}(t_i^J):= f'(t_i^J+0)=\lim _{h\to 0^+}\frac{f(t_i^J+h)-f(t_i^J)}{h}.
$$
Hence, 
$$
\partial f(t_i^J) =\left[\underline{b}(t_i^J), \, \overline{b}(t_i^J) \right].
$$
Let $\{t_i^J\}_{i=1}^k$  be the points of discontinuity between $t_1$ and $t_2$, i.e. 

$ t_2 < t_1^J < \cdots < t_k^J<t_1 $, where $\underline{b}(t_i^J)\leq \overline{b}(t_i^J)$.  
We can estimate as follows

\begin{eqnarray*}
\frac{b(t_1)-b(t_2)}{t_1-t_2} & = & \frac{b(t_1)-b(t_2)+ \displaystyle \sum _{i=1}^k\overline{b}({t^J_i})-\displaystyle \sum _{i=1}^k\overline{b}(t^J_i)}{t_1-t_2} \\
&\geq& \frac{b(t_1)-b(t_2)+ \displaystyle \sum _{i=1}^k\underline{b}(t^J_i)-\displaystyle \sum _{i=1}^k\overline{b}(t^J_i)}{t_1-t_2}\\
&= &\frac{b(t_1)-\overline{b}(t^J_k) -b(t_2)+\underline{b}(t^J_1)
+ \displaystyle \sum _{i=2}^k \big(\underline{b}(t^J_i)-\overline{b}(t^J_{i-1})\big)}{t_1-t_2} \\
&\geq& \frac{-\alpha_0(t_1-t^J_k)-\alpha_0(t^J_1-t_2)-\alpha_0 \displaystyle \sum _{i=2}^k  (t^J_i -t^J_{i-1})}{t_1-t_2}\\ 
&=& \frac{-\alpha_0(t_1-t_2)}{t_1-t_2}=-\alpha_0,
\end{eqnarray*}
where $\alpha_0 := \displaystyle \max_i{c_{L_i}}$ and 
$c_{L_i}$ is the piecewise Lipschitz constant. Hence, the assumption (\ref{eq1}) follows immediately. 
\end{example}

\section{An intermediate HVI % hemivariational
 formulation of the interface problem}
In this section we provide a first equivalent  weak variational formulation  of the interface problem \reff{a1}~-~\reff{a5} 
in terms of a hemivariational inequality (HVI). Since this 
HVI lives on the unbounded domain $\Omega \times \Omega^c$ 
(as the original problem), this HVI cannot numerically 
treated directly and therefore provides only 
an intermediate step in the solution procedure.

For the bounded Lipschitz domain $\Omega$  we use the standard Sobolev space $H^s(\Omega)$ 
and the Sobolev spaces on the bounded Lipschitz boundary $\Gamma$  
(see \cite[Sect 2.4.1]{SaSch2011}),
\begin{eqnarray*} H^s(\Gamma) = 
    \left\{ \begin{array}{ll}
        \{ u|_{\Gamma} : u\in  H^{s+1/2}(\R^d) \} & (0 < s  \le 1) ,  \\
         L^2(\Gamma) & (s = 0),  \\
        (H^{-s}(\Gamma))^* \text{ (dual space) } & (-1 \le s <0).
    \end{array} \right.
  \end{eqnarray*}
Further we need for the unbounded domain $\Omega^c=\R^d \backslash \overline{\Omega}$ the Frechet space
 (see e.g. \cite[Section 4.1, (4.1.43)]{HsWe-2008})
$$ H^s_{loc}(\Omega^c)  = \{u \in {\cal D }^*(\Omega^c) :
 \chi u \in H^s(\R^d) \,\, \forall \chi \in C_0^\infty (\Omega^c) \} \,.
$$ 

By the trace theorem we have $u|_{\Gamma}\in H^{1/2}(\Gamma)$ 
for $u \in H^1_{loc}(\Omega^c)$.
Next we define
$\Phi : H^1(\Omega)\times H^1_{loc}(\Omega^c) \to \R \cup\{\infty\} $
by
\begin{equation}\label{b1}
\Phi(u_1,u_2) := \int_{\Omega} g (|\nabla u_1|)\, dx +
                 \frac 12 \int_{\Omega^c} |\nabla u_2|^2\, dx
                - L(u_1,u_2|_\Gamma ).
\end{equation}

Here the data $f_0 \in L^2(\Omega), t_0 \in L^2(\Gamma)$
enter the linear functional
\begin{equation}\label{b2}
L(u,v) := \int_{\Omega} f_0 \cdot u\, dx +
                  \int_{\Gamma} t_0 \cdot v  \, ds\,.
\end{equation}
Further in \reff{b1} the function $g$ is given by $p$ (see \reff{a1}) through
\begin{displaymath}
g:[0,\infty)\to [0,\infty), t\mapsto g(t) = \int_0^t s\cdot p(s)\, ds ,
\end{displaymath}
where we assume that $p$ is $C^1$, $ 0\le p(t) \le p_0<\infty $,
and  $ t \mapsto t\cdot p(t)$ is 
 strictly monotonic increasing.
Then, $ 0 \le g(t)\le \frac{1}{2} p_0 \cdot t^2 $ and 
the real-valued functional 
\begin{displaymath}
G(u):=  \int_{\Omega} g (|\nabla u|)dx,  \qquad  u\in H^1(\Omega)
\end{displaymath}
is strictly convex. The Frechet derivative of $G$,
\[
DG(u;v)=\int_{\Omega} p (|\nabla u|)  (\nabla u)^T \cdot \nabla v \, dx
\qquad  u,v\in H^1(\Omega)
\]
is strongly monotone in $H^1(\Omega)$ with respect to the semi-norm 
$$| v|_{H^1(\Omega)}= \|\nabla \, v \|_{L^2(\Omega)} \,, $$
that is, there exists a constant $c_G>0$ such that 
\begin{equation} \label{DG-coerc}
c_G\: \hnorm{u- v}{H^1(\Omega)}^2  \le DG(u; u-v)-DG(v; u-v)
\quad \forall  u, v \in H^1(\Omega).
\end{equation}
Analogously to \cite{CCJG-1997, MaiSte-2005}  we first define 
\begin{eqnarray*}
 \left. \begin{array}{ll}
\mathcal{L}_0:= \{ v\in H^1_{loc} (\Omega^c) : 
& \Delta v= 0 \; \mbox{in} \; H^{-1} (\Omega^c) \\[0.5ex]
& \mbox{(and for } d=2 \, \exists a \in \R
\; \mbox{such that } v \mbox{ satisfies }  (\ref{a3}))\},
 \end{array} \right. 
\end{eqnarray*}
and then the affine, hence convex set of admissible functions
\[
C := \{
(u_1,u_2)\in H^1(\Omega)\times H^1_{loc}(\Omega^c) \, : \, \;
         u_1|_{\Gamma_t} = u_2|_{\Gamma_t}+ u_0|_{\Gamma_t}  
\mbox{ and }  u_2 \in \mathcal{L}_0 \}.
\]
According to \cite[Remark 4]{CCJG-1997}, $C$ is closed in $H^1(\Omega)\times H^1_{loc}(\Omega^c)$. 
Further, we have 
\begin{eqnarray*}
 \left. \begin{array}{ll}
D\Phi ((\hat{u}_1,\hat{u}_2);(u_1,u_2))= & D G (\hat{u}_1;u_1)\, 
+ \int _{\Omega^c} \nabla \hat{u}_2 \cdot \nabla u_2 \, dx \\[0.5ex]
& - \int _{\Omega} f_0 \cdot u_1 \, dx 
\, - \int _{\Gamma_t} t_0 \cdot u_2|_{\Gamma_t} \, ds \,.
 \end{array} \right. 
\end{eqnarray*}

Now we can pose the HVI problem $(P_\Phi)$: Find $(\hat{u}_1, \hat{u}_2) \in C$ such that for all $(u_1, u_2) \in C$ there holds
for $\delta u_1 := u_1-\hat{u}_1 , \delta u_2 := u_2-\hat{u}_2$,  
\begin{equation} \label{HVI-1}
D \Phi ((\hat{u}_1, \hat{u}_2);( \delta u_1, \delta u_2))
 + J^0 (\gamma (\hat{u}_2- \hat{u}_1 + u_0); 
\gamma (\delta u_2 - \delta u_1)) \geq 0 \,.
\end{equation}

\begin{theorem} \label{th1}
The HVI problem $(P_\Phi)$ is equivalent to (\ref{a1})~-~(\ref{a5})
in the sense of distributions.
\end{theorem}
{\bf Proof.} First, taking into account the definition of the generalized gradient, we note that 
\begin{equation} \label{eq4}
 \int_{\Gamma_s} p(|\nabla \hat{u}_1|)\frac{\partial \hat{u}_1}{\partial n}
~ \psi \, ds \leq 
 \int _{\Gamma_s} j^0 (\cdot, (\hat{u}_2 - \hat{u}_1 + u_0)|_{\Gamma_s};
 \psi) \, ds, \,
 \forall \psi \in C^ \infty(\Gamma_s)
\end{equation}
 is the integral formulation of the nonmonotone boundary inclusion  \reff{a5}.
\\
 Let $(\hat{u}_1,\hat{u}_2)\in C$ solve (\ref{HVI-1}). To show that $(\hat{u}_1,\hat{u}_2)$ solves (\ref{a1})~-~(\ref{a5}) in the sense of distributions, first choose $\eta \in C_0^{\infty}(\R^d)$ such that  
$(u_1, u_2) := (\hat{u}_1 + \eta|_{\Omega}, \hat{u}_2+\eta|_{\Omega^c}) \in C$. Setting $(u_1, u_2)$ in (\ref{HVI-1}) and integration by parts, implies
\begin{eqnarray*}
& 0 \,\leq & -\int_\Omega\Bigl( f_0 +\div p(|\nabla \hat{u}_1|) \nabla \hat{u}_1 \Bigr) \cdot \eta \, dx
  -  \int_{\Omega^c} \Delta  \hat{u}_2\cdot \eta \, dx \\[0.75ex]
& & + \ds
\int_{\Gamma}\Bigl(p(|\nabla \hat{u}_1|)\frac{\partial \hat{u}_1}{\partial n}
   -\frac{\partial \hat{u}_2}{\partial n}-t_0 \Bigr) \cdot  \eta \, ds \,, 
  \end{eqnarray*}
since the last term in (\ref{HVI-1}) vanishes for the chosen $(u_1,u_2)$, and moreover, $n$ pointing into $\Gamma_0 $ yields the negative sign of
$\frac{ \partial \hat{u}_2}{\partial n} $.  
Varying $\pm \eta\in C_0^{\infty} (\Omega) $ and
$\pm \eta\in C_0^{\infty} (\Omega^c) $,
shows that \reff{a1} and \reff{a2} hold in the sense of distributions. Hence,
\begin{displaymath}
0\le \int_{\Gamma}
\Bigl(  p(|\nabla \hat{u}_1|)\frac{\partial \hat{u}_1}{\partial n}
 - \frac{\partial \hat{u}_2}{\partial n}-t_0  \Bigr) \cdot  \eta \, ds .
\end{displaymath}
Since $\eta $ is arbitrary on $\Gamma$, we obtain 
\begin{equation} \label{eq4a}
p(|\nabla \hat{u}_1|)\frac{\partial \hat{u}_1}{\partial n} =
\frac{\partial \hat{u}_2}{\partial n}+t_0 \qquad \mbox{a.e. on} \; \Gamma .
\end{equation}
This proves (\ref{a5_1}) and the second relation in (\ref{a4}).

Next, let $\eta_1,\eta_2\in C_0^{\infty}(\R^d) $ and consider
$(u_1, u_2):=(\hat{u}_1+ \eta_1|_\Omega, \hat{u}_2+ \eta_2|_{\Omega^c})\in C$ in an analogous way to obtain
\[
0 \leq \int_\Gamma  p(|\nabla \hat{u}_1|)\frac{\partial \hat{u}_1}{\partial n} \eta_1 -(\frac{\partial \hat{u}_2}{\partial n}+t_0) \eta_2 \, ds 
+ \int _{\Gamma_s} j^0 (\cdot,(\hat{u}_2- \hat{u}_1 + u_0)|_{\Gamma_s}; 
(\eta_2-\eta_1)|_{\Gamma_s}) \, ds,
\]
which implies by \reff{eq4a}
\[
 \int_{\Gamma} p(|\nabla \hat{u}|)\frac{\partial \hat{u}_1}{\partial n}
 ~ (\eta_2-\eta_1) \, ds \leq 
 \int _{\Gamma_s} j^0 (\cdot,(\hat{u}_2 - \hat{u}_1 + u_0)|_{\Gamma_s}; 
(\eta_2-\eta_1)|_{\Gamma_s}) \, ds \,. 
\]
Finally, we define $\psi:=\eta_1-\eta_2$. Taking $\eta_1=\eta_2$ on $\Gamma_t$, we have $\psi=0$ on $\Gamma_t$, but $\psi$ is arbitrary on $\Gamma_s$, what
gives (\ref{eq4}). 

Vice versa we show that (\ref{HVI-1}) follows from (\ref{a1})~-~(\ref{a5}). Let $(u_1, u_2)$ solve (\ref{a1})~-~(\ref{a5}). Due to (\ref{a2}) and (\ref{a3}), $(u_1, u_2) \in C$. Multiplying (\ref{a1}) and (\ref{a2}) with 
differences $w_1:= v_1 - u_1$, $w_2:= v_2 - u_2$, respectively, 
where one chooses $(v_1,v_2) \in C$ arbitrarily,
and integrating by parts yields
\begin{eqnarray} 
DG(u_1;w_1) -\int_{\Gamma} p(|\nabla u_1|) 
\frac{\partial u_1}{\partial n} \, w_1 |_{\Gamma}  ds 
& = & \int_{\Omega} f_0 \cdot w_1 \, dx, \label{eq3_1} \\[0.5ex]
\int_{\Omega^c} \nabla u_2 \cdot \nabla w_2 \, dx + \int_{\Gamma} 
\frac{\partial u_2}{\partial n} \, w_2 |_{\Gamma} ds 
& = &  0 \,. \label{eq3_2}
\end{eqnarray}
Combining (\ref{eq3_1}), (\ref{eq3_2}), 
and $p(|\nabla u_1|) \frac{\partial u_1}{\partial n} - \frac{\partial u_2}{\partial n} = t_0$ on $\Gamma$,
we obtain
\[
D\Phi ((u_1,u_2); (w_1, w_2))=\int_{\Gamma} p(|\nabla u_1|) 
\frac{\partial u_1}{\partial n} \, (w_1 - w_2)|_{\Gamma} \, ds \,,
\]
where the latter integral vanishes on $\Gamma_t$ by definition of $C$.
Hence, by (\ref{a5}) (see in particular the integral  
formulation (\ref{eq4})) we conclude that for all $(v_1,v_2) \in C$,
$w_1 = v_1 - u_1$, $w_2 = v_2 - u_2$,  
\begin{eqnarray*}
& & D \Phi ((u_1, {u}_2);(w_1, w_2))
+ \int_{\Gamma_s} j^0 (\cdot,(u_2 - u_1 + u_0)|_{\Gamma_s};
(w_2 - w_1)|_{\Gamma_s})~ds \\
  &  & =  - \int_{\Gamma_s} p(|\nabla u_1|) 
\left.\frac{\partial u_1}{\partial n} (w_2 - w_1)) \right \vert_{\Gamma_s} \, ds 
\\ & & \; \; \; + \int_{\Gamma_s} j^0 (\cdot,(u_2 - u_1 + u_0)|_{\Gamma_s};
(w_2 - w_1)|_{\Gamma_s})~ds
 \\ & & \geq 0 \,,
\end{eqnarray*}
what shows that $(u_1, u_2)\in C$ solves (\ref{HVI-1}). \qed

\begin{remark}
A solution of the HVI problem $(P_\Phi)$ is a critical point of the potential 
energy function 
$$
\Pi(u_1, u_2 ):=  \Phi(u_1, u_2)  + J(\gamma(u_2 - u_1 + u_0))\,,  
$$
which is nonsmooth and nonconvex. Thus generally the 
following  nonsmooth, nonconvex constrained optimization problem: 
\begin{eqnarray}
\label{energy}
\begin{array}{ll}
\mbox{minimize} & \Pi(u_1, u_2) \\
\mbox{subject to} & (u_1, u_2) \in C, \,
\end{array}
\end{eqnarray}

admits several local minimizers which are critical points of $\Pi$
in the sense of  (\ref{HVI-1}).
\end{remark} 

% Under the  smallness condition  \reff{cc2}
% together with  \reff{cc1}, by Corollary \ref{cor-1}, 
% the HVI problem $(P_\Phi)$
%and so also the optimization problem (\ref{energy}) is uniquely solvable.
% $\Omega^c=\R^d \backslash \overline{\Omega}$ NUR Frechet Raum!
% Corollary \ref{cor-1} nur anwendbar auf $\Omega^c \cap K(0;R)$ 
% und dann Fortsetzung wie bei ODEs mit Abstrahlbedingung? 

\section{The boundary/domain HVI 
formulation of the interface problem}

In this section  we employ boundary integral operator theory  \cite{HsWe-2008,GwiSte-2018} 
to rewrite  the exterior problem 
(\ref{a2})~-~(\ref{a3}) as a boundary variational inequality on $\Gamma$. As a result we arrive at an equivalent hemivariational formulation of the original interface problem \reff{a1}~-~\reff{a5} that lives on 
 $\Omega \times \Gamma$ and consists of a weak formulation of the nonlinear differential operator in the bounded domain $\Omega$, the 
Poincare-St\'{e}klov operator on the bounded boundary $\Gamma$, and a nonsmooth functional on the boundary part $\Gamma_s$.

To this end we need the following representation formula, see 
\cite[(1.4.5)]{HsWe-2008},\cite[(12.28)]{GwiSte-2018}.
\begin{lemma}
\label{l1}
For $u_2 \in{\cal L}_0 $ with Cauchy data $(v,\psi)$  $ $
there holds
\begin{equation}\label{c1}
u_2(x') = \frac 12 ( K v(x') - V \psi(x') ) + a, \,
x' \in\Omega^c \,, 
\end{equation}
where $V$ and $K$ denote the single layer potential and the 
double layer potential, respectively,
and $a$ is the constant appearing in \reff{a3} for $d=2$ (and
$a=0$ in (\ref{c1}) if $d > 2 $). 
\end{lemma}
Note that \reff{c1} determines $u_2$ in $\Omega^c $ as far as one knows
its Cauchy data on $\Gamma$.

Next we recall the Poincar\'{e}--Steklov operator for the exterior problem, $S:  H^{1/2}(\Gamma) \rightarrow H^{-1/2}(\Gamma) $ 
is a selfadjoint operator with the defining property
$$
 \partial_n u_2|_\Gamma = - S(u_2|_\Gamma)
$$
for solutions $u_2 \in \mathcal{L}_0 $  of the Laplace equation on $\Omega^c$.
The operator $S$ enjoys the important property that it can be expressed as 
\[
S = \frac{1}{2} [ W  +
 (I-K' ) V^{-1}  (I-K) ] \,,
\]
where $I, V, K, K', W$ denote the identity, the single layer boundary integral operator, the double layer boundary integral operator, its formal adjoint, and the hypersingular integral operator, respectively.
Further, $S$ gives rise to the positive definite bilinear form
$\langle S \cdot, \cdot \rangle$, 
that is, there exists a constant $c_S >0$ such that 
\begin{equation}\label{pos-def}
\langle S v, v \rangle \geq  c_S \|v\|^2_{H^{1/2}(\Gamma)} \,,
\quad  \forall v \in H^{1/2}(\Gamma) \,,
\end{equation}
where $\langle \cdot, \cdot \rangle$ extends the $L^2$ duality on $\Gamma$.

Let $E:= H^1(\Omega)\times \widetilde H^{1/2}(\Gamma_s )$
with $\widetilde H^{1/2}(\Gamma_s ) := \{ w\in H^{1/2}(\Gamma) |
\textnormal{supp}\, w\subseteq \overline\Gamma_s \} $.

We mainly follow \cite{MaiSte-2005,GMS-2021} and use the affine change of variables
\begin{equation}\label{trafo}
(u_1,u_2) \mapsto (u,v) := (u_1 - c, u_0 + u_2|_\Gamma - u_1|_\Gamma) \in E 
\end{equation}
for a suitable $c \in \R \, (c=0,\mbox{if }d=2)$. Note that $v$ is indeed supported in $ \overline\Gamma_s$, since the boundary condition
$(u_1 - u_2)|_{\Gamma_t} = u_0$ guarantees $v|_{\Gamma_t} = 0$.

Next, define the linear functional $\lambda \in E^*$  by 
\begin{displaymath}
\lambda(u,v) := 
 \int _{\Omega} f_0 \cdot u \, dx
 + \langle t_0 + Su_0, u|_{\Gamma} + v\rangle \,, \quad (u,v) \in E 
\end{displaymath}

and consider the potential energy function 
$$
{\cal P}(u,v) := G(u) + 
\frac{1}{2} \langle S(u|_\Gamma +v), u|_\Gamma  + v \rangle 
+ J(v) - \lambda(u,v) 
= \Pi(u_1, u_2 ) + C \,,
%:=  \Phi(u_1, u_2)  + J(j(u_2 - u_1 + u_0))  
$$
where $ C = C(u_0,t_0)$ is a constant independent of $u,v$.

Moreover, in case $d=2$, similarly to \cite{CCJG-1997}, we 
introduce an additional linear constraint 
and consider the affine closed subspace
\begin{displaymath}
D := \{ (u,v)\in E \, : \,  
\langle S1, u|_\Gamma + v - u_0 \rangle = 0 \; \mbox{if }d=2  \; \}. 
\end{displaymath}

Then the nonsmooth, nonconvex constrained optimization problem: 
\begin{eqnarray}
\label{Energy}
\begin{array}{ll}
\mbox{minimize} & {\cal P}(u, v) \\
\mbox{subject to} & (u, v) \in D \,
\end{array}
\end{eqnarray}
leads to the following hemivariational inequality problem $(P_\mathcal{A})$:
 Find   $(\hat u,\hat v)\in D $ such that for all $(u,v) \in D$, 
\begin{equation} \label{HVIPJ}
\mathcal{A}(\hat{u}, \hat{v}; u-\hat{u}, v-\hat{v}) 
+ J^0 (\hat{v}; v-\hat{v}) \, ds \geq \lambda (u-\hat{u}, v-\hat{v}) \,, 
\end{equation}
where $\mathcal{A} : E \to E ^*$  is defined for all $(u,v), (u',v') \in E$ by 
\begin{displaymath}
\mathcal{A}(u,v)\, (u',v') = \mathcal{A}(u,v; u',v')
:= DG(u,u') + \langle S(u|_\Gamma + v),  u'|_{\Gamma} + v' \rangle \,.
\end{displaymath}

Since the nonfrictional convex smooth part of $\cal P$ coincides with the functional $J$ of \cite{MaiSte-2005}, we immediately obtain from
\cite[Theorem 2]{MaiSte-2005} that the  problems $(P_\Phi)$ and $(P_\mathcal{A})$ are equivalent in the following way. 

\begin{theorem}\label{p2}
(i) Let $(u_1,u_2)\in C $ solve $(P_\Phi)$.
Then, $(u,v)\in D$ defined by (\ref{trafo}) solves $(P_\mathcal{A})$.
\\ 
(ii) Let  $(u,v)\in D$ solve $(P_\mathcal{A})$.
Take $a\in\R$ arbitrarily if $d=2$, whereas $a=0$ if $d = 3$. 
Define $u_1:=u + a$ and $u_2$ by the representation formula (\ref{c1})
with $(u|_{\Gamma}+v-u_0+a,-S(u|_{\Gamma}+v-u_0))$ replacing $(v,\psi)$, i.e. 
$u_2:=\frac{1}{2} \left( K(u|_{\Gamma} + v-u_0+a)
+ V(S(u|_{\Gamma}-v+u_0 ))\right) + a$.
Then, $(u_1,u_2)\in C$ solves $(P_\Phi)$.
\end{theorem}

Thanks to the strong monotonicity of the nonlinear operator $DG$  in $H^1(\Omega)$ with respect to the semi-norm $| \cdot|_{H^1(\Omega)}= \|\nabla \, \cdot \|_{L^2(\Omega)}$ and the positive definiteness of the 
Poincar\'{e}--Steklov operator $S$ the following strong monotonicity 
property can be derived, see \cite[Lemma 4.1]{CCJG-1997}. 

\begin{lemma} \label{l-coerc}
There exists a constant $c_0>0$
such that for all
$v,v' \in \tilde H^{1/2}(\Gamma_s)$ and all $u,u' \in H^1(\Omega)$
there holds
\begin{eqnarray*}
\lefteqn{ c_0 \cdot \Vert (u-u',v-v')
\Vert^2_{ H^1(\Omega)\times \tilde H^{1/2}(\Gamma_s) }}\\
&\le &  DG(u; u-u')- DG(u';u-u') \\ &+&
 \langle S(u|_\Gamma + v - u'|_\Gamma - v')}{u|,
u_\Gamma + v - u'|_\Gamma - v' \rangle.
\end{eqnarray*}
\end{lemma}

Moreover, following the arguments in \cite{OvGw-Rassias-14} and
using the compact embedding 
$\widetilde H^{1/2}(\Gamma_s ) \subset L^2 (\Gamma_s)$
 it can be easily seen that the  functional 
$\varphi : \widetilde H^{1/2}(\Gamma_s ) \times \widetilde H^{1/2}(\Gamma_s )
 \to \R$ defined by
\[
 \varphi (v,\tilde{v}) 
:= \int _{\Gamma_s} j^0(\cdot, v;  \tilde{v} - v)\, ds
\]
is pseudo-monotone and upper semicontinuous.   
Thus the concrete hemivariational inequality problem $(P_\mathcal{A})$
is covered by the general theory exhibited in Section \ref{inter};
existence and under the smallness condition uniqueness hold for $(P_\mathcal{A})$.

\section{A regularized  boundary/domain HVI % hemivariational
formulation of the interface problem}
\label{reg}
In this section, we  first recall from \cite{OvGw-2014, OvGw-Rassias-14} a smoothing approximation for a class of nonsmooth functions that can be expressed by means of the  plus function $\hat p(x)=x^+=\max \{x, 0 \}$. Then we apply this smoothing approximation to the nonsmooth locally Lipschitz function $j(s,\cdot)$ and thus arrive at a regularized formulation of the hemivariational inequality problem $(P_\mathcal{A})$, for which we can provide existence and uniqueness results. Finally in this section, a convergence result for the developed regularization procedure is presented.  

The general idea of a smoothing approximation is to use convolution.  
However in general, convolution  is not easily applicable in practice.
On the other hand, for a special class of functions that can be expressed by means of the plus function, a smoothing approximation can be explicitly computed as follows.  Let us first focus to the case $f(x)=\max\{g_1(x),g_2(x)\}$, where  $g_1, g_2$ are smooth functions.  Then,  
\begin{equation} \label{pr1}
f(x) = g_1(x) + \hat p[ g_2(x)- g_1(x)].
\end{equation}
With the notations 
\[
\R_{+}=\{\varepsilon\in \R \, :\, \varepsilon \geq 0\}, \quad \R_{++}=\{\varepsilon\in \R \, :\, \varepsilon > 0\}
\]
we define the smoothing approximation $P \, :\, \R_{++} \times \R$  of $\hat p$ via convolution by
\[
P(\varepsilon, x)=\int_{\R} {\hat p}(x-\varepsilon t) \rho(t) \, dt.
\] 
Here, $\varepsilon >0$ is a small regularization parameter and 
$\rho \, :\, \R \to \R_{+}$ is a probability density function such that
$$
\kappa= \int_{\R} |t|\, \rho(t)\, dt < \infty.
$$
  Now we replace $\hat p(x)$ by its approximation $P(\varepsilon, x)$ and obtain  $\hat{f}:\R_{++}\times \R \to \R$, 
\begin{equation} \label{smooth}
\hat{f}(\varepsilon,x) =g_1(x)+ P(\varepsilon,g_2(x)-g_1(x)),
\end{equation}
as a smoothing function of $f$ in (\ref{pr1}). 

 Using, for example, the Zang probability density function \cite{Zang}
\[
\rho(t)= \left \{ \begin{array} {ll} 1 & \mbox{if} \, - \frac{1}{2}\leq t \leq
    \frac{1}{2} \\
0 & \mbox{otherwise},
\end{array} \right. 
\]
leads to  
\begin{equation} \label{smooth_Zang}
P(\varepsilon, x) = \int _{\R} \hat p(x-\varepsilon t) \rho(t) \, dt = \left \{ \begin{array} {ll} 0 & \mbox{if} \quad x < -\frac {\varepsilon}{2}\\ 
\frac{1}{2\varepsilon}(x+ \frac{\varepsilon}{2})^2 & \mbox{if} \, -
\frac{\varepsilon}{2} \leq x \leq \frac{\varepsilon}{2} \\
x  &\mbox{if} \quad x > \frac
    {\varepsilon}{2}
\end{array}
\right. 
\end{equation}
and hence, 
\[
\hat{f}(\varepsilon, x): = 
\left \{ \begin{array} {ll} g_1(x) & \mbox{if} \, \;
    (i) \; \mbox{holds}  %\nonumber
     \\[0.1ex]
\frac{1}{2\varepsilon} [g_2(x) -g_1(x)]^2  + 
\frac{1}{2} [g_2(x) + g_1(x)] + \frac{\varepsilon}{8} & \mbox{if} \, \; (ii)  \; \mbox{holds}  
\\ %\label{smoothZang}
g_2(x) & \mbox{if} \, \; (iii) \;  \mbox{holds};
\end{array}
\right.
\]
where the cases 
$(i)$,  $(ii)$,  $(iii)$ are defined,  respectively, by
\begin {description}
\item (i) $g_2(x) -g_1(x) \leq - \frac{\varepsilon}{2} $
\item (ii) $- \frac{\varepsilon}{2} \leq
g_2(x) -g_1(x) \leq \frac{\varepsilon}{2}$
\item (iii) $g_2(x) -g_1(x) \geq \frac{\varepsilon}{2}$.
\end{description}
For other choices of probability density functions and other examples of smoothing functions we refer to \cite{OvGw-Rassias-14, OvGw-2014} and the references therein.

We can extend the representation formula (\ref{pr1}) to the maximum function $f:\R \to \R$ of $m$ smooth functions $g_1, \ldots, g_m$. Indeed we write 
\begin{eqnarray*} %\label{eq:bsp_f_function}
f(x)&=&\max \{g_1(x), g_2(x), \ldots, g_m(x)\} \\
&=& g_1(x) + \hat p[g_2(x)-g_1(x)+ \ldots + \hat p[g_m(x)-g_{m-1}(x)]].
\end{eqnarray*}
The smoothing function $\hat{f} :\R_{++} \times \R \to R$ is then given by
\begin{equation} \label{S_general}
\hat{f}(\varepsilon, x)= 
g_1(x) + P\left( \varepsilon, g_2(x)-g_1(x)+\ldots + 
P\left(\varepsilon, 
g_m(x)-g_{m-1}(x) \right)\right).
\end{equation}
More general nested max-min functions can be treated 
by an appropriate representation using the plus function, 
what is exploited in the numerical treatment of related
nonsmmoth variational problems, albeit in bounded domains, see 
\cite{GwOv-2020}. 
Note that this smoothing procedure via the plus function needs only one regularization parameter $\varepsilon$.

The major properties of the function $\hat{f}(\cdot, \cdot)$ in (\ref{S_general}) 
are listed in the following lemma:  
\begin{lemma} \cite{OvGw-2014,OvGw-Rassias-14}\label{major-reg}
\begin {description}
\item (i) For any $\varepsilon >0$ and for all $x\in \R$, 
\[
|\hat{f}(\varepsilon,x)-f(x)|\leq (m-1)\kappa \varepsilon.
\]
\item(ii) The function $\hat{f}$ is continuously differentiable on $\R_{++}\times \R$ and for any $x\in \R$ and $\varepsilon>0$ there exist $\Lambda_i \in [0,1]$ such that $\displaystyle \sum_{i=1}^m \Lambda_i=1$ and 
 \begin {equation}\label{f1}
 \frac{\partial \hat{f}(\varepsilon,x)}{\partial x} = \hat{f}_x( \varepsilon, x)=\displaystyle \sum_{i=1}^m \Lambda_i  g'_i(x).
\end{equation}
Moreover,
\begin {equation}\label{f2}
\{\limsup_{z\to x, \varepsilon \to 0^+}  \hat{f}_x( \varepsilon, z)\}\subseteq \partial f(x).
\end{equation}
\end {description}
\end{lemma}

Now analogously let the superpotential 
 $j \, : \, \Gamma_s \times \R \to \R $ be a maximum function of smooth functions, i.e.
$$
j(s,x) = \max \{g_1(s,x), g_2(s,x), \ldots, g_m(s,x)\}
$$
and obtain the  smoothing function $\hat{j} \, :\, \Gamma_s \times
\R_{++} \times \R \to R$ via
\begin{equation} \label{S_j_general}
\hat{j}(s,\varepsilon, x)= 
g_1(s,x) + P\left( \varepsilon, g_2(s,x)-g_1(s,x)+\ldots + 
P\left(\varepsilon, 
g_m(s,x)-g_{m-1}(s,x) \right)\right).
\end{equation}

Assume that for $i= 1, \ldots, m$
there exist positive constants $c_i, d_i$ such that for
almost all (a.a) $s\in \Gamma_s$ and for all $x\in \R$,
\begin{equation} \label{ass_12}
| g'_i(s,x)| \leq  c_i(1+|x|) , \,  g'_i(s,x) \, x \geq  -d_i|x|. 
\end{equation}
Under (\ref{ass_12}), the growth condition  \reff{as-j}
 in Section \ref{inter} is immediately satisfied. Moreover, from Lemma \ref{major-reg} %(\ref{f1})-(\ref{f2})
 and (\ref{ass_12})  the following auxiliary result  can be easily deduced. 
\begin{lemma}  \label{auxiliary}
It holds for all  $x, z, \xi \in \R$, almost all (a.a.)
$s \in \Gamma_s$  that
\begin{eqnarray} 
&& \left| \hat{j}_x (s, \varepsilon, x) \, z \right | \leq  c(1+|x|)\,|z|  \label{ass_3}\\
 &&  \hat{j}_x (s, x, \varepsilon) \cdot (-x)  \leq  d|x| \label{ass_4}\\
&&  \limsup_{z\to x, \varepsilon\to 0^+} \hat{j}_x(s, \varepsilon, z)\, \xi 
\leq  j^0(s, x; \xi) \,.  \label{ass_3b}
\end{eqnarray}
\end{lemma}
Next we introduce  $J_\varepsilon :  H^{1/2}(\Gamma_s)  \to \R$  by
$$
J_\varepsilon(v)=\int_{\Gamma_s} \hat{j}(s, \varepsilon, v(s)) \, ds.
$$
Since $\hat{j}(s,\varepsilon, \cdot) $ is continuously differentiable, the functional $J_\varepsilon$ is G\^{a}teaux differentiable with continuous G\^{a}teaux derivative 
$D J_\varepsilon :  H^{1/2}(\Gamma_s) \to  (H^{1/2}(\Gamma_s) )^*$,
\[
\langle D J_\varepsilon({v}), \tilde{v} \rangle 
= \int_{\Gamma_s} \hat{j}_x (s, \varepsilon, v(s)) ~ 
\tilde{v} (s) \, ds.
\]

The {\it regularized problem} ($P_\varepsilon$) of ($P_{\mathcal{A}}$) reads now: 
Find $({\hat u}_{\varepsilon}, \hat{v}_{\varepsilon}) \in D \subset
E = H^1(\Omega) \times \tilde H^{1/2}(\Gamma_s)$ 
such that 
\begin{equation} \label{reg-2D} 
\mathcal{A} (\hat{u}_{\varepsilon}, \hat{v}_{\varepsilon})
(u -\hat{u}_{\varepsilon}, v-\hat{v}_{\varepsilon}) + 
\langle D J_\varepsilon({\hat v}_{\varepsilon}), v - \hat{v}_{\varepsilon} \rangle
 =  \lambda (u-\hat{u}_{\varepsilon}, v -\hat{v}_{\varepsilon})
 \quad \forall (u,v) \in D.
\end{equation}
Note that $D$ is an affine subset of 
$E = H^1(\Omega) \times \tilde H^{1/2}(\Gamma_s)$  
and $J_\varepsilon$ is differentiable; therefore
 ($P_\varepsilon$) simplifies to the variational equality (\ref{reg-2D}).
%, which in 3D simplifies further to the operator equation

Assume there exists a constant $0\leq \alpha <c_S$, where $c_S$ is defined in (\ref{pos-def}), such that for a.a. $s \in \Gamma_s$ and for all
$x_1, \,x_2 \in \R,$
that
\begin{equation}\label{uniq_reg}
(\hat{j}_x(s, \varepsilon, x_1)-\hat{j}_x(s,\varepsilon, x_2))(x_1-x_2)
\geq -\alpha |x_1-x_2|^2 \,, 
\end{equation}
then by the discussion in Section \ref{inter}
the regularized problem ($P_\varepsilon$) has a unique solution 
$\mathbf{u}_\varepsilon  = (u_{\varepsilon}, v_{\varepsilon})$.

Moreover replacing the coercive bilinear form $a$ in \cite{OvGw-2014},
by the strongly monotone form, see Lemma \ref{l-coerc},
\begin{eqnarray*}
&& [(u,v), (u',v')] \in E \times E
\mapsto \\
&&  DG(u; u-u')- DG(u';u-u') \\
&& +
 \langle S(u|_\Gamma + v - u'|_\Gamma - v'),
u_\Gamma + v - u'|_\Gamma - v' \rangle \,,
\end{eqnarray*}
the approximation result \cite[Theorem 5.1]{OvGw-2014} for the regularization procedure extends to the present situation. Thus we can summarize our findings as follows.   

\begin{theorem} \label{theo:exist+str-conv} 
The regularized problem ($P_\varepsilon$) has at least one solution $\mathbf{u}_\varepsilon \in D \subset
 E = H^1(\Omega) \times \tilde H^{1/2}(\Gamma_s) $. The family $\{\mathbf{u}_\varepsilon\}$ is uniformly bounded in $E$.
 Moreover, there exists a subnet of $\{\mathbf{u}_\varepsilon\}$ which converges strongly in $E$ to a solution $\mathbf{u}$ of the hemivariational inequality  problem $(P_\mathcal{A})$. 
Under the smallness condition $0\leq \alpha <c_S$
along with the condition (\ref{uniq_reg}) strong convergence of the entire regularization procedure $\{\mathbf{u}_\varepsilon\}$ follows.
\end{theorem}

\section{Numerical approximation} \label{num-approx}
In this section we treat the numerical approximation for
problem $(P_\varepsilon)$ by a Galerkin
projection using finite elements in $\Omega $ and boundary
elements on $\Gamma $ leading to the coupling of FEM and BEM.
 Here we focus to the case $d=2$ and assume the interior domain $\Omega$ to be polygonal.  For the more general case of a domain with a piecewise $C^{1,1}$ boundary we can refer to the construction of an intermediate polygonal 
approximation with finite element/boundary element
coupling and the associated error analysis presented in
\cite[Sections 3.3 - 3.5]{gatih95} for a class of nonlinear problems. 
 
Let $(H_h^1\times \widetilde H_h^{1/2} \times  H_h^{-1/2})_{h \in I}$
be a family of finite dimensional subspaces of
$H^1(\Omega)\times \widetilde H^{1/2}(\Gamma_s)\times
H^{-1/2}(\Gamma) $ where $I \subseteq (0,\infty)$ with $0\in \bar I $.
Let us specify the finite dimensional ansatz spaces $H_h^1$,  $ H_h^{-1/2}$ and $\widetilde H_h^{1/2}$ as follows.
With the  mesh parameter $h>0$, we
consider a nested regular quasi-uniform family $({\cal T}_h)_h$ 
of meshes of $\Omega$ consisting of triangles, say, with a diameter between $c_1 \cdot h $ 
and $c_2 \cdot h~(0<c_1<c_2) $,
and denote by $H_h^1$ the space of the continuous and piecewise linear trial functions associated to the triangulation ${\cal T}_h$, that is,
\[
H_h^1:= \{u_h \in C(\bar{\Omega}) \, : \,u_h |_{T} \in \GP_1 (T) \quad \forall \, T \in \mathcal{T}_h\} \subset H^1(\Omega).
\]
The family $({\cal T}_h)_h$ of meshes of $\Omega $ leads to the mesh family 
$({\cal E}_h)_h$ of edges on the boundary $\Gamma$ so that we may introduce
$H_h^{1/2}$ and 
$H_h^{-1/2} $ as the  space of all continuous and piecewise linear functions,
respectively, as the space of all piecewise constant functions on $\Gamma$, associated to the partition ${\cal E}_h$.
Further, we assume that the partition  ${\cal E}_h$ of the boundary leads to
a partition $\widetilde{\mathcal{E}}_h$ of $\Gamma_s$, so that
$\widetilde H^{1/2}_h $  becomes the subspace of those trial functions in $H_h^{1/2}$ that are supported on $\Gamma_s$.
Thus we have
\[
{H}^{-1/2}_h := \{ w_h \in {L}^2(\Gamma) \, : \, w_h|_{E}  \in \GP_0 (E) \quad \forall E \in {\mathcal {E}}_h\} \subset {H}^{-1/2}(\Gamma)
\]
and 
\[
\widetilde{H}^{1/2}_h :=\{{v_h} \in {C}(\Gamma_s)\, :\, {v_h}|_{E} \in \GP_1 (E) \quad \forall E \in \widetilde{\mathcal{E}}_h\} 
\subset \widetilde{H}^{1/2}(\Gamma_s).
\]
Hence, we may simply define  
$E_h:= H_h^1\times  \widetilde{H}_h^{1/2}$.

For $h\in I $ let $i_h : H_h^1\hookrightarrow H^1(\Omega),
j_h : H_h^{1/2} \hookrightarrow  H^{1/2}(\Gamma),$
and $k_h : H_h^{-1/2} \hookrightarrow H^{-1/2}(\Gamma)$
% IF NEEDED:
% $ \widetilde j_h := j_h|\widetilde H_h^{1/2} :
% \widetilde H_h^{1/2}  \hookrightarrow \widetilde H^{1/2}(\Gamma_s) $
denote the canonical imbeddings with their duals
$i_h^*$, $j_h^*$, and $k_h^*$.

A straightforward discretization of the problem  $(P_\varepsilon)$
would lead to replace
%$E$ by $E_h=H_h^1\times \widetilde H_h^{-1/2} $  and hence
$S$ by $ j_h^* S j_h $. 
Unfortunately, the numerical calculation of
$j_h^*  S j_h$ requires
the numerical computation of $V^{-1}$. Since $V^{-1}$ is, in general, not known
explicitly we have to approximate $V^{-1} $ and therefore introduce the
subspace $H^{-1/2}_h $ to approximate the tractions on the interface. 
Here we follow the Costabel symmetric BEM/FEM coupling procedure and approximate $S$ by 
\begin{displaymath}
S_{h} := \frac{1}{2} [
j_h^*  W j_h +
j_h^*  (I-K' )k_h (k_h^* V k_h)^{-1} k_h^* (I-K) j_h ] \,.
\end{displaymath}
The computation of $ S_{h} $ requires the numerical
solution of a linear system with a symmetric, positive definite
matrix $ V_{h}:=(k_h^* V k_h)$.
In general  $ S_h \ne j_h^*  S  j_h$, 
since $S_h$ is the Schur complement of a discretized matrix, what only needs
the knowledge of $(k^*_h V k_h)^{-1}$,  while
$j_h^* S j_h$ is the  discretized Schur complement of operators,
what needs the knowledge of $ V^{-1}$. 

Here we recall two essential properties of $ S_h $.

\begin{lemma} \cite[Lemma 12.2.9]{GwiSte-2018}
\label{cl2}
There exist a constant $c_0>0$ and $h_0\in I $ such that for any
$h\in I  $ with $h< h_0 $ and any $ v_h \in H_h^{1/2}$  there holds
\begin{equation}\label{d10}
\xdba{S_h v_h}{v_h} \:\ge \: c_0 \:
\norm{ v_h }{ H^{1/2}(\Gamma) }^2. 
\end{equation}
\end{lemma}

\begin{lemma} \cite[Lemma 5.2]{CCJG-1997},\cite[(3.15)]{CoSt-1990}
\label{cl3}
There exists a constant $c>0$ such that for all $h\in I $ and
for all $v_h\in H_h^{1/2}$ there holds
\begin{displaymath}
\norm{S_h v_h - j_h^* S j_h v_h }{H^{-1/2}(\Gamma)} \le
c \: \dist_{H^{-1/2}(\Gamma) }\Bigl(  V^{-1} (I-K) j_h v_h ;
 H_h^{-1/2}\Bigr).
\end{displaymath}
\end{lemma}

Let $\mathcal{D}_h$ be another partition of $\Gamma_s$  consisting of elements $K_i$ joining the midpoints $P_{i-1/2}$, $P_{i+1/2}$ of the edges $E \in \widetilde{\mathcal{E}}_h$ lying on $\Gamma_s$ sharing $P_{i}$ as a common point. If $P_i$ is a vertex of $\partial \Omega$ then $K_i$ is the half of the edge. On $\mathcal{D}_h$ we introduce the space $\mathcal{Y}_h$ of all piecewise constant functions by 
\[
\mathcal{Y}_h=\{\mu_h \in L^{\infty}(\Gamma_s)\, :\, \mu_h|_K \in \GP_0(K) \quad \forall K \in \mathcal{D}_h\}
\]
and define the piecewise constant Lagrange interpolation operator $
L_h \,:\, \widetilde{H}^{1/2}_h \to \mathcal{Y}_h$ by
\[
L_h(w_h)(x)=\sum _{i} %{P_i \in \Gamma_C \cap \Sigma_h}
 w_h(P_i) \,
\chi_{\mbox{\small int} \, _{\Gamma_s} K_i}(x),
\]
where $\chi_{\small \mbox{int} \,_{\Gamma_s} K_i}$ is the characteristic function of the interior of $K_i$ in $\Gamma_s$.
Then we employ the approximation 
\begin{equation} \label{exint}
\langle D J_{\varepsilon,h}(v_h), w_h \rangle 
= \displaystyle  {\int_{\Gamma_s}}  
\hat{j}_x(s, \varepsilon, L_h(v_h)(s)) ~ L_h(w_h)(s)  \,ds. 
\end{equation}
Further we note that there holds with some positive constant $c$ not dependent on $h$, for any $v_h \in H_h^{1/2}$ 
\begin{eqnarray} \label{gl-1}
&& |L_h v_h\|_{L^2(\Gamma)} \leq c\|v_h\|_{L^2(\Gamma)} \\
\label{gl-2}
&& \|L_h v_h -v_h\|_{L^2(\Gamma)} \leq ch^{1/2} \|v_h\|_{H^{1/2}(\Gamma)},  
\end{eqnarray}
what follows by real interpolation from standard estimates, see e.g. \cite{Gl-2008}.

To simplify matters, we assume that the datum $u_0$ 
in the transmission condition (\ref{a4}) belongs to
$H_h^{1/2}$ for all $0<h<h_0$ and so 
$D_h := E_h \cap D$.

Now we can state the discrete version $(P_{\varepsilon,h})$ associated to the 
regularized problem $(P_{\varepsilon})$:
Find $(\hat u_{\varepsilon,h}, \hat v_{\varepsilon,h})  \in D_h $ such that 
\begin{eqnarray} \label{disc-reg-2D} 
&& {\cal A}_h (\hat u_{\varepsilon,h},\hat v_{\varepsilon,h} ;u_{h}- \hat u_{\varepsilon,h},v_{h}- \hat v_{\varepsilon,h}) 
  + \,  \langle DJ_{\varepsilon,h} (\hat{v}_{\varepsilon,h}), v_h - \hat{v}_{\varepsilon,h} \rangle \\ \nonumber 
&& = % \ge
\lambda_h(u_h- \hat u_{\varepsilon,h},v_h- \hat v_{\varepsilon,h})
\quad \forall (u_h,v_h)\in D_h, 
\end{eqnarray}
where
${\cal A}_h:H^1_h\times\widetilde H^{1/2}_h\to(H^1_h\times\widetilde H^{1/2}_h)^*$
is defined via
\begin{eqnarray*}
&& {\cal A}_h (u_h,v_h) (r_h,s_h) = {\cal A}_h (u_h,v_h; r_h,s_h) \\
 && := DG (u_h; r_h) +
\langle S_h (u_h|_\Gamma + v_h), r_h|_\Gamma + s_h \rangle 
\end{eqnarray*}
and $\lambda_h : H_h^1 \times \widetilde{H}_h^{1/2} \to \R$  by
\begin{eqnarray*}
&  \lambda_h (u_h, v_h) : =   & \int_\Omega f_0 \cdot u_h \, dx 
 + \langle t_0, u_h |_{\Gamma} - v_h \rangle \\
&& + \langle i_h^* \gamma^*(W+(I-K')k_h(k_h^*V k_h)^{-1}k_h^*(I-K))u_0,
u_h|_{\Gamma} - v_h \rangle.
\end{eqnarray*}

The solvability of $(P_{\varepsilon,h})$ follows from the existence result in \cite[Theorem 3.1]{GwOv-2015} for general pseudo-monotone bifunctions.
Similar to Theorem \ref{theo:exist+str-conv} we can guarantee uniqueness
of solutions $\mathbf{\hat u}_{\varepsilon,h} := 
(\hat u_{\varepsilon,h}, \hat v_{\varepsilon,h}) $
to $(P_{\varepsilon,h})$. Moreover, since $DG$ is 
strongly monotone in $H^1(\Omega)$ with respect to the semi-norm
 $|\cdot|_{H^1(\Omega)}= \|\nabla \, \cdot \|_{L^2(\Omega)}$, 
see  (\ref{DG-coerc}), and $S_h$ is uniformly coercive in $H^{1/2}(\Gamma)$,
see Lemma \ref{cl2}, (\ref{d10}), 
the uniform boundedness of the family 
$(\mathbf{\hat u}_{\varepsilon,h})_h$
can be derived. To sum up, we have the following result.  

\begin{theorem} \label{th-exist+unibound}
Under the smallness condition $0\leq \alpha <c_S$
along with the condition (\ref{uniq_reg}), 
for $h<h_0$ ($h_0 $ given in Lemma \ref{cl2})
the problem $(P_{\varepsilon, h})$ has  exactly one solution.
Moreover, the family $(\mathbf{\hat u}_{\varepsilon,h})_h$
is uniformly bounded in $E$. 
\end{theorem}

Finally we can prove the following a priori error estimate.
\begin{theorem} \label{th-apriori}
For $(\hat{u}_\varepsilon, \hat{v}_\varepsilon)\in D$, the solution of problem $(P_\varepsilon)$, and $(\hat{u}_{\varepsilon,h},\hat{v}_{\varepsilon,h}) \in D_h$, the solution of problem $(P_{\varepsilon, h})$, there holds 
for $h<h_0$  and  $\alpha> 0$ sufficiently small (in (\ref{uniq_reg}))
with a constant $C>0$, which is
 independent of $h$, but depends on $\alpha$,  
%, but depending on $\varepsilon>0$
\begin{align}  %{eqnarray}
 & C \, \|(\hat{u}_\varepsilon-\hat{u}_{\varepsilon,h},
\hat{v}_\varepsilon -
\hat{v}_{\varepsilon,h})\|^2_{H^1(\Omega)\times \widetilde{H}^{1/2}(\Gamma_s)} 
 \label{apriori} \\  
% &\leq    \inf_{(u_h,v_h) \in D_h}
% \; \left\{
% \|\hat{u}_\varepsilon-u_{h}\|^2_{H^1(\Omega)} 
% + \|\hat{v}_\varepsilon-v_h\|^2_{\tilde{H}^{1/2}(\Gamma_s)} \right. \nonumber
% \\
% &  \left. \quad  + \,\, \|\hat{v}_\varepsilon-v_h\|_{\tilde{H}^{1/2}(\Gamma_s)}  \right \} \nonumber \\
&\leq    \inf_{(u_h,v_h) \in D_h}
 \; \left\{
\|\hat{u}_\varepsilon-u_{h}\|^2_{H^1(\Omega)} 
+ \|\hat{v}_\varepsilon-v_h\|^2_{\tilde{H}^{1/2}(\Gamma_s)} 
  +   \|\hat{v}_\varepsilon-v_h\|_{L^2(\Gamma)}  
 \right \} 
 \nonumber \\    
 & +  \mbox{dist}_{H^{-1/2}(\Gamma)}
 (V^{-1}(I-K)(\hat{u}_\varepsilon+\hat{v}_\varepsilon-u_0), H_h^{-1/2})^2  
 + h^{1/2}  \nonumber \,.           
\end{align} %{eqnarray}
\end{theorem}
{\bf Proof.} 
First, similar to the proof of Theorem 5 in \cite{MaiSte-2005}, 
in virtue of Lemma \ref{cl2}, (\ref{d10}), 
further by  (\ref{disc-reg-2D}) and (\ref{reg-2D}), 
there holds for $0<h<h_0$ 
modulo  a positive constant, independent of $h$, 
for all $(u_h, v_h) \in D_h$ 
\begin{eqnarray} 
&& \| u_h|_\Gamma - \hat{u}_{\varepsilon,h}|_\Gamma +
v_h -\hat{v}_{\varepsilon,h})\|^2_{H^{1/2}(\Gamma)} \nonumber\\  
&\lesssim & 
\langle S_h (u_h|_\Gamma - \hat{u}_{\varepsilon,h}|_\Gamma +
v_h -\hat{v}_{\varepsilon,h}),
u_h|_\Gamma - \hat{u}_{\varepsilon,h}|_\Gamma +
v_h -\hat{v}_{\varepsilon,h} \rangle \nonumber \\
&=&
\langle S_h (u_h|_\Gamma - \hat{u}_{\varepsilon,h}|_\Gamma +
v_h -\hat{v}_{\varepsilon,h}),
u_h|\Gamma - \hat{u}_{\varepsilon,h}|_\Gamma +
v_h -\hat{v}_{\varepsilon,h} \rangle \nonumber \\
&& +
 ({\cal A}_h (\hat u_{\varepsilon,h},\hat v_{\varepsilon,h}) - \lambda_h)
 (u_{h}- \hat u_{\varepsilon,h},v_{h}- \hat v_{\varepsilon,h}) \nonumber \\
 && +
 ({\cal A} (\hat u_{\varepsilon},\hat v_{\varepsilon}) - \lambda)
 (\hat u_{\varepsilon,h} - \hat u_{\varepsilon},  
\hat v_{\varepsilon,h} - \hat v_{\varepsilon} ) \nonumber \\
   && + \,\, 
\langle DJ_{\varepsilon,h}(\hat{v}_{\varepsilon,h}), 
v_h - \hat{v}_{\varepsilon,h}\rangle + 
\langle DJ_\varepsilon(\hat{v}_\varepsilon), 
\hat v_{\varepsilon,h} - \hat v_{\varepsilon} \rangle
\,. \label{apriori-1}
\end{eqnarray}
To conclude the a priori error estimate (\ref{apriori}),
inspecting the proof of Theorem 5 in \cite{MaiSte-2005},
 we only need to estimate the last both terms in (\ref{apriori-1}). 
For this purpose  we use the shorthand 
$ \hat{f}_x(v) :=  \hat{j}_x(\varepsilon, \cdot, v)$
and proceed as follows:
\begin{eqnarray}
 && 
\langle DJ_{\varepsilon,h}(\hat{v}_{\varepsilon,h}),v_h - \hat{v}_{\varepsilon,h} \rangle
+ \langle DJ_\varepsilon(\hat{v}_\varepsilon),
\hat v_{\varepsilon,h} - \hat v_{\varepsilon} \rangle \nonumber  \\ 
&=&
\int_{\Gamma_s} \hat{f}_x(L_h \hat{v}_{\varepsilon,h})
L_h  (v_h  - \hat{v}_{\varepsilon,h} ) \, ds + 
\int_{\Gamma_s} \hat{f}_x(\hat{v}_\varepsilon) \,
(\hat v_{\varepsilon,h} - \hat v_{\varepsilon})
   \, ds \nonumber   \\ 
&=& \int_{\Gamma_s} \hat{f}_x(L_h \hat{v}_{\varepsilon,h}) \,
 ( L_h v_h  - \hat{v}_\varepsilon) \, ds
  + \int_{\Gamma_s} \hat{f}_x (\hat{v}_{\varepsilon}) \,
 (\hat{v}_{\varepsilon,h} - L_h \hat{v}_{\varepsilon,h}) \, ds \nonumber \\
&& +  \int_{\Gamma_s} \left(
\hat{f}_x(L_h \hat{v}_{\varepsilon,h}) - \hat{f}_x(\hat{v}_{\varepsilon}) \right) \, \left( \hat{v}_{\varepsilon} -  L_h \hat{v}_{\varepsilon,h} \right) \, ds \,.
 \label{apriori-2}
\end{eqnarray}
Next we analyze the three summands in (\ref{apriori-2}) separately.
For the first, we obtain by Lemma \ref{auxiliary}, (\ref{ass_3}),
\[
\int_{\Gamma_s} \hat{f}_x(L_h \hat{v}_{\varepsilon,h}) \,
 ( L_h v_h  - \hat{v}_\varepsilon) \, ds
\lesssim   \| 1 + | L_h \hat{v}_{\varepsilon,h}| \|_{L^2(\Gamma_s)} 
 \| L_h v_h  - \hat{v}_\varepsilon \|_{L^2(\Gamma_s)} \,,
\]
further for the first factor, using the estimate (\ref{gl-1}),
\begin{align*}
 \|1 + | L_h \hat{v}_{\varepsilon,h}| \|_{L^2(\Gamma_s)}
\lesssim (1 +  \| \hat{v}_{\varepsilon,h} \|_{L^2(\Gamma_s)}),
\end{align*}
what is uniformly bounded by Theorem \ref{th-exist+unibound},
and for the second factor,
by triangle inequality and by (\ref{gl-2}),
\begin{align*} 
& \| L_h v_h  -\hat{v}_\varepsilon \|_{L^2(\Gamma)}
\le \| L_h v_h  - v_h  \|_{L^2(\Gamma)} 
+ \| v_h  -\hat{v}_\varepsilon \|_{L^2(\Gamma)} \\
&  \lesssim h^{1/2} \| v_h \|_{H^{1/2}(\Gamma)} 
+ \| v_h  -\hat{v}_\varepsilon \|_{L^2(\Gamma)} \,. 
\end{align*} 
For the second summand, using the estimate (\ref{gl-2}),
\begin{align*}
\| \hat{v}_{\varepsilon,h} - L_h \hat{v}_{\varepsilon,h} \|_{L^2(\Gamma_s)}
 \lesssim h^{1/2} \, \|\hat{v}_{\varepsilon,h}\|_{H^{1/2}(\Gamma_s)}, 
\end{align*}
where $\|\hat{v}_{\varepsilon,h}\|_{H^{1/2}(\Gamma_s)}$  is uniformly bounded by Theorem \ref{th-exist+unibound}.
For the third summand we use the assumption (\ref{uniq_reg}) and estimate
\begin{align*}
 \int_{\Gamma_s} \left(
\hat{f}_x(L_h \hat{v}_{\varepsilon,h}) - \hat{f}_x(\hat{v}_{\varepsilon}) \right) \, \left( \hat{v}_{\varepsilon} -  L_h \hat{v}_{\varepsilon,h} \right) \, ds \le \alpha \|  \hat{v}_{\varepsilon} -  L_h \hat{v}_{\varepsilon,h} \|^2_{L^2(\Gamma_s)} \,,
\end{align*}
then by triangle inequality and further by (\ref{gl-2})
and Theorem \ref{th-exist+unibound}, just as above,
\begin{align*}
& \|  \hat{v}_{\varepsilon} -  L_h \hat{v}_{\varepsilon,h} \|^2_{L^2(\Gamma_s)} 
\le 2 \| \hat{v}_{\varepsilon,h} -  L_h \hat{v}_{\varepsilon,h} \|^2_{L^2(\Gamma_s)}
+  2 \|  \hat{v}_{\varepsilon} -  \hat{v}_{\varepsilon,h}  \|^2_{L^2(\Gamma_s)}  \\
& \lesssim 
( h + \|\hat{v}_{\varepsilon} -  \hat{v}_{\varepsilon,h}\|^2_{H^{1/2}(\Gamma_s)} )\,.
\end{align*}

Finally, we combine the above estimates together with (\ref{apriori-1}) and (\ref{apriori-2}) and thus can extend the error estimate
 of Theorem 5 in \cite{MaiSte-2005} to 
arrive at (\ref{apriori}) for small enough $\alpha$. 
\qed

We remark that  Theorem \ref{th-apriori} extends the error estimates in \cite[Theorem 4]{CCJG-1997} and \cite[Theorem 5]{MaiSte-2005} to the nonlinear interface problem with the nonmonotone transmission condition, including the $h$-approximation of the regularized functional $DJ_\varepsilon$ by $DJ_{\varepsilon,h}$. 

\section{Conclusions and an Outlook} 

This paper has shown how various techniques from different fields of mathematical analyis and numerical analysis can be combined to a solution procedure for a nonlinear interface problem that models nonmonotone frictional
contact of elastic infinite media.

Here we have considered the coercive situation. A more delicate problem arises with a loss of coercivity, when only unilateral friction conditions, but no classical transmission conditions in the form of equalities are prescribed on the coupling boundary. For such semicoercive/noncoercive problems, albeit on a bounded domain, we can refer to  e.g. \cite{GwOv-2020}.

The errror analysis in this paper has focused to the classical $h$-method of 
BEM/FEM discretization, for the more versatile $hp$-BEM and $hp$-FEM
applied to unilateral contact problems, albeit in bounded domains,  
we refer to \cite{Ov-Banz} and to \cite{Gwi-13}, respectively.   

Another direction of research is the coupling of the BEM with other discretizations methods for the interior problem, like discontinuous 
Galerkin methods and mixed finite element methods; see e.g.
\cite{GatMaiSte-2011} for the numerical analysis of a transmission problem with Signorini contact. 

\bibliographystyle{amsplain}

\bibliography{biblio-coupling-BE-FE-smoothing}

\end{document}